\documentclass[a4paper,10pt]{article}
\usepackage[frenchb,english]{babel}
\usepackage{amssymb}
\usepackage{amsmath}

\newcommand{\HH}{{\mathbb H}}
\newcommand{\CC}{{\mathbb C}}
\newcommand{\RR}{{\mathbb R}}
\newcommand{\PP}{{\mathbb P}}

\newcommand{\SSS}{{\mathbb S}}

\newcommand{\cZ}{{\mathcal Z}}
\renewcommand{\phi}{\varphi}

\newcommand{\half}{\frac{1}{2}}

\newcommand{\Spin}{\mathrm{Spin}}

\renewcommand{\Re}{\mathrm{Re\ }}
\newcommand{\End}{\mathrm{End}}

\newcommand{\tr}{\mathrm{tr}}

\newcommand{\vol}{\mathrm{vol}}

\newcommand{\SZ}{{\Sigma \cZ}}
\newcommand{\SZp}{{\Sigma^+ \cZ}}

\newcommand{\Spinc}{\mathrm{Spin^c}}

\newcommand{\Ric}{\mathrm{Ric}}

\newcommand{\DDM}{\widetilde D}

\newcommand{\PSZ}{P_{\SSS^1}\cZ}
\newcommand{\PSM}{P_{\SSS^1}M}

\newcommand{\cercle}{\mathbb{S}}
\newtheorem{thm}{Theorem}[section]

\newtheorem{remark}[thm]{Remark}
\newtheorem{remarks}[thm]{Remarks}
\newtheorem{definition}[thm]{Definition}
\newtheorem{notation}[thm]{Notation}
\newtheorem{example}[thm]{Example}

\begin{document}
\title{{\bf The Energy-Momentum tensor on low dimensional $\Spinc$ manifolds}}

\author{Georges Habib\\
{\small\it Lebanese University,
Faculty of Sciences II, Department of Mathematics}\\ 
{\small \it P.O. Box 90656 Fanar-Matn, Lebanon}\\
{\small \it ghabib@ul.edu.lb}\\\\ 
Roger Nakad\\
{\small \it Max Planck Institute for Mathematics,
Vivatsgasse 7, 53111 Bonn, Germany}\\
{\small \it nakad@mpim-bonn.mpg.de}}

\date{}
\maketitle
\noindent On a compact surface endowed with any $\Spinc$ structure, we give a formula involving the Energy-Momentum tensor in terms of geometric quantities. A new proof of a B\"{a}r-type inequality for the eigenvalues of the Dirac operator is given. The round sphere $\mathbb{S}^2$ with its canonical $\Spinc$ structure satisfies the limiting case. Finally, we give a spinorial characterization of immersed surfaces in $\mathbb{S}^2\times \mathbb{R}$ by solutions of the generalized Killing spinor equation associated with the induced $\Spinc$ structure on $\mathbb{S}^2\times \mathbb{R}$.\\\\
{\it Keywords:} $\Spinc$ structures, Dirac operator, eigenvalues, Energy-Momentum tensor, compact surfaces, isometric immersions.\\ \\
Mathematics subject classifications (2000): $53$C$27$, $53$C$40$, $53$C$80$.

\section{Introduction} 
\setcounter{equation}{0}
On a compact $\Spin$ surface, Th. Friedrich and E.C. Kim proved that any eigenvalue $\lambda$  of the Dirac operator satisfies the equality \cite[Thm. 4.5]{fk}:
\begin{equation}\label{eq:Fk}
\lambda^2=\frac{\pi\chi(M)}{Area(M)}+\frac{1}{Area(M)}\int_M|T^\psi|^2v_g, 
\end{equation}
where  $\chi(M)$ is the Euler-Poincar\'e characteristic of $M$ and $T^\psi$ is the field of quadratic forms called the Energy-Momentum tensor. It is given on the complement set of zeroes of the eigenspinor $\psi$  by 
$$ T^\psi (X, Y) = g(\ell^\psi (X), Y) = \frac 12 \Re(X\cdot\nabla_Y\psi + Y\cdot\nabla_X\psi,\frac{\psi}{\vert\psi\vert^2}),$$
for every $X, Y \in \Gamma(TM)$. Here $\ell^\psi$ is the field of symmetric endomorphisms associated with the field of quadratic forms $T^\psi$. We should point out that since $\psi$ is an eigenspinor, the zero set is discret \cite{barrr}.
The proof of Equality (\ref{eq:Fk}) relies mainly on a local expression of the covariant derivative of $\psi$ and the use of the Schr\"odinger-Lichnerowicz formula. This equality has many direct consequences.
First, since the trace of $\ell^\psi$ is equal to $\lambda$, we have  by the Cauchy-Schwarz inequality that $\vert \ell^\psi\vert^2 \geqslant \frac{(tr(\ell^\psi))^2}{n} = \frac{\lambda^2}{2},$ where $tr$ denotes the trace of $\ell^\psi$. Hence, Equality (\ref{eq:Fk}) implies the B\"{a}r inequality \cite{bar2} given by
\begin{eqnarray}
\lambda^2 \geqslant \lambda_1^2 := \frac{2\pi\chi(M)}{Area(M)}.
 \label{barspin}
\end{eqnarray}
Moreover, from Equality (\ref{eq:Fk}), Th. Friedrich and E.C. Kim \cite{fk} deduced
 that $\int_M {\rm det}(T^\psi)v_g=\pi\chi(M)$, which gives an information on the Energy-Momentum tensor without knowing the eigenspinor nor the eigenvalue. Finally, for any closed surface $M$ in $\mathbb{R}^3$ of constant mean curvature $H$, the restriction to $M$ of a parallel spinor on $\mathbb{R}^3$ is a generalized Killing spinor of eigenvalue $-H$ with Energy-Momentum tensor equal to the Weingarten tensor $II$ (up to the factor $-\frac{1}{2})$ \cite{m1} and we have by \eqref{eq:Fk}
$$H^2=\frac{\pi\chi(M)}{Area(M)}+\frac{1}{4Area(M)}\int_M|II|^2v_g.$$
Indeed, given any surface $M$ carrying such a spinor field, Th. Friedrich \cite{fr3} showed 
that the Energy-Momentum tensor associated with this spinor  satisfies the Gauss-Codazzi equations and hence $M$ is locally immersed into $\mathbb{R}^3$.\\\\
Having a $\Spinc$ structure on manifolds is a weaker condition than having a $\Spin$ structure because  every $\Spin$ manifold has a trivial $\Spinc$ structure. Additionally, any compact surface or any product of a compact surface with $\mathbb{R}$ has a $\Spinc$ structure carrying particular spinors. In the same spirit as in \cite{13}, when using a suitable conformal change, the second author \cite{r1} established a B\"{a}r-type inequality for the eigenvalues of the Dirac operator on a compact surface endowed with any $\Spinc$ structure. In fact, any eigenvalue $\lambda$ of the Dirac operator satisfies
\begin{eqnarray}
\lambda^2 \geqslant \lambda_1^2 := \frac{2\pi \chi(M)}{Area(M)} - \frac{1}{Area(M)} \int_M \vert\Omega\vert v_g,
\label{barspinc}
\end{eqnarray}
where $i\Omega$ is the curvature form of the connection on the line bundle given by the $\Spinc$ structure.
Equality is achieved if and only if the eigenspinor $\psi$ associated with the first eigenvalue $\lambda_1$ is a Killing $\Spinc$ spinor, i.e., for every $X\in \Gamma(TM)$ the eigenspinor $\psi$ satisfies
\begin{eqnarray}
\left\{
\begin{array}{l}
\nabla_X\psi = -\frac{\lambda_1}{2}X\cdot\psi,\\
\Omega\cdot\psi =i \vert\Omega\vert\psi.
\end{array}
\label{lamba}
\right.
\end{eqnarray}
Here $X\cdot\psi$ denotes the Clifford multiplication and $\nabla$ the spinorial Levi-Civita connection \cite{fr1}.\\\\
Studying the Energy-Momentum tensor on a compact Riemannian $\Spin$ or $\Spinc$
manifolds has been done by many authors, since it is related to several geometric situations.
Indeed, on compact $\Spin$ manifolds, J.P. Bourguignon and P. Gauduchon \cite{BG92}
proved that the Energy-Momentum tensor appears naturally in the study of the variations
of the spectrum of the Dirac operator. Th. Friedrich and E.C. Kim \cite{fk2} obtained
the Einstein-Dirac equation as the Euler-Lagrange equation of a certain functional. The second author extented 
these last two results to $\Spinc$ manifolds \cite{r2}. Even if it is not a
computable geometric invariant, the Energy-Momentum tensor is, up to a constant, the
second fundamental form of an isometric immersion into a $\Spin$ or $\Spinc$ manifold carrying a
parallel spinor \cite{m1, r2}. For a better understanding of the tensor $q^\phi$
 associated with a spinor field $\phi$, the first author \cite{habib1}
 studied Riemannian flows and proved
that, if the normal bundle carries a parallel spinor $\psi$, the tensor $q^\phi$ associated with $\phi$ (the restriction of $\psi$ to the flow) 
 is the O'Neill tensor of the flow.\\\\
In this paper, we give a formula corresponding to \eqref{eq:Fk} for any eigenspinor $\psi$ of the square of the Dirac operator on compact surfaces endowed with any $\Spinc$ structure (see Theorem \ref{procarrdir}). It is motivated by the following two facts: First, when we consider eigenvalues of the square of the Dirac operator, another tensor field is of interest. It is the skew-symmetric tensor field  $Q^\psi$ given by 
$$Q^\psi(X,Y)= g(q^\psi (X), Y) = \frac{1}{2}\Re(X\cdot\nabla_Y\psi-Y\cdot\nabla_X\psi,\frac{\psi}{|\psi|^2}),$$ 
for all vector fields $X,Y\in \Gamma(TM).$ This tensor was studied by the first author in the context of Riemannian flows \cite{habib1}. Second, we consider any compact surface $M$ immersed in $\cercle^2\times \mathbb{R}$ where $\cercle^2$ is the round sphere equipped with a metric of curvature one. The $\Spinc$ structure on $\cercle^2\times \mathbb{R},$ induced from the canonical one on $\cercle^2$ and the $\Spin$ struture on $\mathbb{R},$ admits a parallel spinor \cite{mo}. The restriction to $M$ of this $\Spinc$ structure is still a $\Spinc$ structure with a generalized Killing spinor \cite{r2}.\\\\
In Section \ref{pre}, we recall some basic facts on $\Spinc$ structures and the restrictions of these structures to hypersurfaces. In Section \ref{sect:3} and after giving a formula corresponding to \eqref{eq:Fk} for any eigenspinor $\psi$ of the square of the Dirac operator, we deduce a formula for the integral of the determinant of $T^\psi+ Q^\psi$ and we establish a new proof of the B\"{a}r-type inequality (\ref{barspinc}). In Section \ref{sec4}, we consider the $3$-dimensional case and treat examples of hypersurfaces in $\CC \PP^2$. In the last section, we come back to the question of a spinorial characterisation of  surfaces in $\cercle^2\times \mathbb{R}$. Here we use a different approach than the one in \cite{r}. In fact, we prove that given any surface $M$ carrying a generalized Killing spinor associated with a particular $\Spinc$ structure, the Energy-Momentum tensor satisfies the four compatibility equations stated by B. Daniel \cite{daniel}. Thus there exists a local immersion of $M$ into $\cercle^2\times \mathbb{R}$.
\section{Preliminaries}\label{pre}
\setcounter{equation}{0}
In this section, we begin with some preliminaries concerning $\Spinc$ structures and
the Dirac operator. Details can be found in \cite{lm}, \cite{montiel}, \cite{fr1}, \cite{r1} and
\cite{r2}.\\ \\
{\bf The Dirac operator on $\Spinc$ manifolds:} Let $(M^n, g)$ be a Riemannian manifold of dimension $n\geqslant2$ without
boundary. We denote by ${\rm SO}M$ the 
${\rm SO}_n$-principal bundle over $M$ of positively oriented orthonormal frames. A
$\Spinc$ structure of $M$ is a $\Spin_n^c$-principal bundle $(\Spinc M,\pi,M)$
 and an $\cercle^1$-principal bundle $(\cercle^1 M ,\pi,M)$ together with a double
covering given by  $\theta: \Spinc M \longrightarrow {\rm SO}M\times_{M}\cercle^1 M$ such
that $\theta (ua) = \theta (u)\xi(a),$
for every $u \in \Spinc M$ and $a \in \Spin_n^c$, where $\xi$ is the $2$-fold
covering of $\Spin_n ^c$ over ${\rm SO}_n\times \cercle^1$. 
Let $\Sigma M := \Spinc M \times_{\rho_n} \Sigma_n$ be the associated spinor bundle
where $\Sigma_n = \CC^{2^{[\frac n2]}}$ and $\rho_n : \Spin_n^c
\longrightarrow  \End(\Sigma_{n})$ denotes the complex spinor representation. A section of
$\Sigma M$ will be called a spinor field. The spinor bundle $\Sigma M$ is equipped with a
natural Hermitian scalar product denoted by $(.,.)$. We
define an $L^2$-scalar product
$<\psi,\phi> = \int_M (\psi,\phi) v_g,$
for any spinors $\psi$ and $\phi$. \\ Additionally, any connection 1-form $A: T(\cercle^1 M)\longrightarrow i\RR$ on
$\cercle^1 M$ and the connection 1-form 
$\omega^M$ on ${\rm SO} M$, induce a connection
on the principal bundle ${\rm SO} M\times_{M} \cercle^1 M$, and hence 
a covariant derivative $\nabla$ on $\Gamma(\Sigma M)$ \cite{fr1,r2}. The curvature
of $A$ is an imaginary valued 2-form denoted by $F_A= dA$, i.e., $F_A = i\Omega$,
where $\Omega$ is a real valued 2-form on $\cercle^1 M$. We know
 that $\Omega$ can be viewed as a real valued 2-form on $M$ $\cite{fr1, kn}$. In
this case $i\Omega$ is the curvature form of the associated line bundle $L$. It is
the complex line bundle associated with the $\cercle^1$-principal bundle via the
standard representation of the unit circle.
For every spinor $\psi$, the Dirac operator is locally defined 
by   $$D\psi =\sum_{i=1}^n e_i \cdot \nabla_{e_i} \psi,$$
where $(e_1,\ldots,e_n)$ is a local oriented orthonormal tangent frame and ``$\cdot$''  denotes the Clifford multiplication. The Dirac
operator is an elliptic, self-adjoint operator with respect to the $L^2$-scalar
product and verifies, for any spinor field $\psi$, the Schr\"odinger-Lichnerowicz formula 
\begin{eqnarray}
D^2\psi=\nabla^*\nabla\psi+\frac{1}{4}S\psi+\frac{i}{2}\Omega\cdot\psi
\label{bochner}
\end{eqnarray}
where $\Omega\cdot$ is the extension of the Clifford multiplication to differential
forms given by 
$(e_i ^* \wedge e_j ^*)\cdot\psi = e_i\cdot e_j \cdot\psi$. For any spinor $\psi \in \Gamma(\Sigma M)$, we have \cite{mo2}
\begin{eqnarray}
 (i\Omega\cdot\psi,\psi)\ \geqslant -\frac{c_n}{2} \vert \Omega \vert_g\vert
\psi\vert^2,
\label{cs}
\end{eqnarray}
where $\vert \Omega \vert_g$ is the norm of  $\Omega$, with respect to $g$ given by
$\vert \Omega \vert_g^2=\sum_{i<j} (\Omega_{ij} )^2$
in any orthonormal local frame and $c_n = 2[\frac
n2]^\frac12$. Moreover, equality holds in (\ref{cs}) if and only
if $\Omega\cdot\psi = i \frac {c_n}{2}\vert\Omega\vert_g \psi$.\\
Every $\Spin$ manifold has a trivial $\Spinc$ structure \cite{fr1, bm}. In fact, we
choose the trivial line bundle with the trivial connection whose curvature $i\Omega$
is zero. Also every K\"ahler manifold $M$ of complex dimension $m$ has a canonical $\Spinc$ structure. Let $\ltimes$ by the K\"{a}hler form defined by the complex structure $J$, i.e. $\ltimes (X, Y)= g(JX, Y)$ for all vector fields $X,Y\in \Gamma(TM).$ The complexified cotangent bundle 
$$T^*M\otimes \CC = \Lambda^{1,0} M \oplus \Lambda^{0,1}M$$
decomposes into the $\pm i$-eigenbundles of the complex linear extension of the complex structure. Thus, the spinor bundle of the canonical $\Spinc$ structure is given by $$\Sigma M = \Lambda^{0,*} M =\oplus_{r=0}^m \Lambda^{0,r}M,$$
where $\Lambda^{0,r}M = \Lambda^r(\Lambda^{0,1}M)$ is the bundle of $r$-forms of type $(0, 1)$. The line bundle of this canonical $\Spinc$ structure is given by  $L = (K_M)^{-1}= \Lambda^m (\Lambda^{0,1}M)$, where $K_M$ is the canonical bundle of $M$ \cite{fr1, bm}. This line bundle $L$ has a canonical holomorphic connection induced from the Levi-Civita connection whose curvature form is given by $i\Omega = -i\rho$, where $\rho$ is the Ricci form given by $\rho(X, Y) = \Ric(JX, Y)$. 
We point out that the canonical $\Spinc$ structure on every K\"{a}hler manifold carries  a parallel spinor \cite{fr1, mo}.\\\\ 
{\bf Spin$^c$ hypersurfaces and the Gauss formula:} Let $\cZ$ be an oriented ($n+1$)-dimensional Riemannian $\Spinc$ manifold and $M
\subset \cZ$ be an oriented hypersurface. The manifold $M$ inherits a $\Spinc$
structure induced from the one on $\cZ$, and we have \cite{r2}
$$ \Sigma M\simeq \left\{
\begin{array}{l}
\Sigma \cZ_{|_M} \ \ \ \ \ \ \text{\ \ \ if\ $n$ is even,} \\\\
 \SZp_{|_M}   \ \text{\ \ \ \ \ \ if\ $n$ is odd.}
\end{array}
\right.
$$
Moreover Clifford multiplication by a vector field $X$, tangent to $M$, is given by 
\begin{eqnarray}
X\bullet\phi = (X\cdot\nu\cdot \psi)_{|_M},
\label{Clifford}
\end{eqnarray}
where $\psi \in  \Gamma(\Sigma \cZ)$ (or $\psi \in \Gamma(\SZp)$ if $n$ is odd),
$\phi$ is the restriction of $\psi$ to $M$, ``$\cdot$'' is the Clifford
multiplication on $\cZ$, ``$\bullet$'' that on $M$ and $\nu$ is the unit normal
vector. The connection 1-form defined on the restricted $\SSS^1$-principal bundle $(\PSM :=
\PSZ_{|_M},\pi,M)$, is given by $A= {A^\cZ}_{|_M} : T(\PSM) = T(\PSZ)_{|_M}
\longrightarrow i\RR.$ Then the curvature 2-form $i\Omega$ on the
$\SSS^1$-principal bundle $\PSM$ is given by $i\Omega= {i\Omega^\cZ}_{|_M}$,
which can be viewed as an imaginary 2-form on $M$ and hence as the curvature form of
the line bundle $L^M$, the restriction of the line bundle $L^\cZ$ to $M$. For every
$\psi \in \Gamma(\SZ)$ ($\psi \in \Gamma(\SZp)$ if $n$ is odd), the real 2-forms
$\Omega$ and $\Omega^\cZ$ are related by \cite{r2}
\begin{eqnarray}
(\Omega^\cZ \cdot\psi)_{|_M} = \Omega\bullet\phi -
(\nu\lrcorner\Omega^\cZ)\bullet\phi.
\label{glucose}
\end{eqnarray}
We denote by $\nabla^{\Sigma \cZ}$ the spinorial Levi-Civita connection on $\Sigma
\cZ$ and by $\nabla$ that on $\Sigma M$. For all $X\in \Gamma(TM)$, we have the spinorial Gauss formula \cite{r2}:
\begin{equation}
(\nabla^{\Sigma \cZ}_X\psi)_{|_M} =  \nabla_X \phi + \half II(X)\bullet\phi,
\label{spingauss}
\end{equation}
where $II$ denotes the Weingarten map of the hypersurface. Moreover, Let $D^\cZ$ and $D^M$ be the Dirac operators on $\cZ$ and $M$, after
denoting by the same symbol any spinor and its restriction to $M$, we have
\begin{equation}
\nu\cdot D^\cZ\phi = \DDM\phi +\frac{n}{2}H\phi -\nabla^{\Sigma\cZ}_{\nu}\phi,
\label{diracgauss}
\end{equation}
where $H = \frac1n \tr(II)$ denotes the mean curvature and $\DDM = D^M$ if $n$ is even and $\DDM=D^M \oplus(-D^M)$ if $n$ is odd.
\section{The $2$-dimensional case}\label{sect:3}
\setcounter{equation}{0}
In this section, we consider compact surfaces endowed with any $\Spinc$ structure. We have 
\begin{thm}\label{procarrdir}
 Let $(M^2, g)$ be a Riemannian manifold and $\psi$ an eigenspinor of the
square of the Dirac operator $D^2$ with eigenvalue $\lambda^2$ associated with any $\Spinc$ structure. Then we have 
$$ \lambda^2 = \frac{S}{4}  +  {\vert{T}^{\psi} \vert^2} + \vert{Q}^{\psi} \vert^2 +\Delta f+|Y|^2-2Y(f)+(\frac i2 \Omega\cdot\psi, \frac{\psi}{\vert\psi\vert^2}), $$ 
where $f$ is the real-valued function defined by $f=\frac{1}{2}{\rm ln}|\psi|^2$
and $Y$ is a vector field on $TM$ given by $g(Y,Z)=
\frac{1}{\vert\psi\vert^2} \Re(D\psi, Z\cdot\psi)$ for any $Z \in \Gamma(TM)$.
\label{thm1}
\end{thm}
{\bf Proof.} Let $\{e_1, e_2\}$ be an orthonormal frame of $TM$. Since the spinor bundle $\Sigma M$ is of real dimension $4$, the set $\{\frac{\psi}{|\psi|}, \frac{e_1\cdot\psi}{|\psi|},\frac{e_2\cdot\psi}{|\psi|},\frac{e_1\cdot e_2\cdot\psi}{|\psi|}\}$ is orthonormal with respect to the real product $\Re(\cdot, \cdot)$. The covariant derivative of $\psi$ can be expressed in this frame as 
\begin{eqnarray} 
 \nabla_X\psi = \delta(X) \psi + \alpha(X)\cdot\psi + \beta(X) e_1\cdot
e_2\cdot\psi,
\label{equa1}
\end{eqnarray}
for all vector fields $X,$ where $\delta$ and $\beta$ are $1$-forms and $\alpha$ is a $(1,1)$-tensor field. Moreover $\beta$, $\delta$ and $\alpha$ are uniquely
determined by the spinor $\psi$. In fact, taking the scalar product of (\ref{equa1}) respectively with $\psi, e_1\cdot\psi, e_2\cdot\psi, e_1\cdot e_2\cdot\psi$, we get $\delta =\frac{d(\vert\psi\vert^2)}{2\vert\psi\vert^2}$ and
$$\alpha (X)=- \ell^\psi(X)+q^\psi(X)\quad\text{and}\quad\beta(X)=\frac{1}{\vert\psi\vert^2} \Re(\nabla_X\psi, e_1\cdot e_2\cdot\psi).$$ 
Using (\ref{bochner}), it follows that
$$\lambda^2 =  \frac {\Delta ( \vert \psi \vert^2 )} { 2  \vert \psi \vert^2 } +
\vert \alpha \vert^2 + \vert \beta\vert^2 + \vert \delta\vert^2 +\frac 14 S + (\frac i2 \Omega\cdot\psi, \frac{\psi}{\vert\psi\vert^2}).$$
Now it remains to compute the term $\vert\beta\vert^2$. We have
\begin{eqnarray*}
\vert\beta\vert^2 &=& \frac{1}{\vert\psi\vert^4}\Re(\nabla_{e_1}\psi, e_1\cdot e_2
\cdot\psi)^2+\frac{1}{\vert\psi\vert^4}\Re(\nabla_{e_2}\psi, e_1\cdot e_2
\cdot\psi)^2\\
&=& \frac{1}{\vert\psi\vert^4} \Re(D\psi - e_2 \cdot\nabla_{e_2}\psi, e_2 \cdot\psi)^2
+ \frac{1}{\vert\psi\vert^4} \Re(D\psi - e_1 \cdot\nabla_{e_1}\psi, e_1 \cdot\psi)^2 \\
&=&  g(Y, e_1)^2 +  g(Y, e_2)^2 + \frac{\vert d(\vert\psi\vert^2)\vert^2}{4
\vert\psi\vert^4} - g(Y, \frac{d(\vert\psi\vert^2)}{\vert\psi\vert^2})\\
&=& |Y|^2-2Y(f) + \frac{\vert d(\vert\psi\vert^2)\vert^2}{4
\vert\psi\vert^4},
\end{eqnarray*}
which gives the result by using the fact that $\Delta f=\frac {\Delta ( \vert \psi \vert^2)} {2\vert \psi \vert^2}+\frac{\vert d(\vert\psi\vert^2)\vert^2}{2
\vert\psi\vert^4}.$
\hfill$\square$
\begin{remark} \label{rempropre} Under the same conditions as Theorem \ref{thm1}, if $\psi$ is an eigenspinor of $D$ with eigenvalue $\lambda$, we get
$$ \lambda^2 = \frac{S}{4}  +  {\vert{T}^{\psi} \vert^2} +\Delta f+(\frac i2 \Omega\cdot\psi, \frac{\psi}{\vert\psi\vert^2}).$$ 
In fact, in this case $Y=0$ and 
\begin{eqnarray}
0 &=& \Re(D\psi,e_1\cdot e_2\cdot\psi)=\Re(e_1\cdot\nabla_{e_1}\psi+e_2\cdot\nabla_{e_2}\psi,e_1\cdot e_2\cdot\psi)\nonumber \\
&=&\Re(-e_2\cdot\nabla_{e_1}\psi+e_1\cdot\nabla_{e_2}\psi,\psi) =2Q^\psi(e_1,e_2)|\psi|^2.
\label{qpsi}
\end{eqnarray} 
This was proven by Friedrich and Kim in \cite{fk} for a $\Spin$ structure on $M$.
\end{remark}
In the following, we will give an estimate for the integral $\displaystyle\int_M {\rm det}(T^\psi+Q^\psi) v_g$ in terms of geometric quantities, which has the advantage that it does not depend on the eigenvalue $\lambda$ nor on the eigenspinor $\psi$. This is a generalization of the result of Friedrich and Kim in \cite{fk} for $\Spin$ structures.
\begin{thm} 
Let $M$ be a compact surface and $\psi$ any eigenspinor of $D^2$ associated with eigenvalue $\lambda^2$. Then we have 
\begin{eqnarray}
\int_M {\rm det}(T^\psi+Q^\psi)v_g \geq \frac{\pi\chi(M)}{2}- \frac 14 \int_M \vert\Omega\vert v_g.
\label{dett}
\end{eqnarray}
Equality in \eqref{dett} holds if and only if either $\Omega$ is zero or has constant sign. 
\label{D2}
\end{thm} 
{\bf Proof.} As in the previous theorem, the spinor $D\psi$ can be expressed in the orthonormal frame of the spinor bundle. Thus the norm of $D\psi$ is equal to
\begin{eqnarray}
|D\psi|^2&=&\frac{1}{|\psi|^2}\Re(D\psi,\psi)^2+\frac{1}{|\psi|^2}\sum_{i=1}^2 \Re(D\psi,e_i\cdot\psi)^2+\frac{1}{|\psi|^2}\Re(D\psi,e_1\cdot e_2\cdot\psi)^2\nonumber\\
&=&({\tr}\,T^\psi)^2|\psi|^2+|Y|^2|\psi|^2+\frac{1}{|\psi|^2}\Re(D\psi,e_1\cdot e_2\cdot\psi)^2,
\label{eq:Dirac}
\end{eqnarray}
where we recall that the trace of $T^\psi$ is equal to $-\frac{1}{|\psi|^2}\Re(D\psi,\psi).$ On the other hand, by (\ref{qpsi}) we have that
$\frac{1}{|\psi|^2}\Re(D\psi,e_1\cdot e_2\cdot\psi)^2=2|Q^\psi|^2|\psi|^2.$ Thus Equation \eqref{eq:Dirac} reduces to 
$$\frac{|D\psi|^2}{|\psi|^2}=({\tr}\, T^\psi)^2+|Y|^2+2|Q^\psi|^2.$$
Now with the use of the equality $\Re(D^2\psi,\psi)=|D\psi|^2-{\rm div}\xi,$ where $\xi$ is the vector field given by $\xi=|\psi|^2Y,$ we get 
\begin{equation}
\lambda^2+\frac{1}{|\psi|^2}{\rm div}\xi=({\tr}\, T^\psi)^2+|Y|^2+2|Q^\psi|^2.
\label{eq:Dirac2}
\end{equation} 
An easy computation leads to $\frac{1}{|\psi|^2}{\rm div}\xi={\rm div} Y+2Y(f)$ where we recall that $f=\frac{1}{2}{\rm ln}(|\psi|^2).$ Hence substituting this formula into \eqref{eq:Dirac2} and using Theorem \ref{thm1} yields 
$$\frac{S}{4} + (\frac i2 \Omega\cdot\psi, \frac{\psi}{\vert\psi\vert^2}) + \Delta f+{\rm div} Y= ({\tr} T^\psi)^2+|Q^\psi|^2-|T^\psi|^2=2{\rm det}(T^\psi+Q^\psi).$$
Finally integrating over $M$ and using the Gauss-Bonnet formula, we deduce the required result with the help of Equation (\ref{cs}). Equality holds if and only if $\Omega\cdot\psi = i\vert\Omega\vert \psi$. In the orthonormal frame $\{e_1, e_2\}$, the $2$-form $\Omega$ can be written $\Omega= \Omega_{12} \ e_1 \wedge e_2$, where $\Omega_{12}$ is a function defined on $M$. Using the decomposition of $\psi$ into positive and negative spinors  $\psi^+$ and $\psi^-$, we find that the equality is attained if and only if
$$\Omega_{12} \ e_1 \cdot e_2 \cdot\psi^+ + \Omega_{12}\  e_1 \cdot e_2 \cdot\psi^- = i \vert\Omega_{12} \vert \psi^+ +i \vert\Omega_{12} \vert \psi^-,$$
which is equivalent to say that,
$$\Omega_{12} \psi^+ = -\vert\Omega_{12}\vert \psi^+\ \ \ \ \text{and}\ \ \ \  \Omega_{12} \psi^- = \vert\Omega_{12}\vert \psi^-.$$
Now if $\psi^+ \neq 0$ and $\psi^- \neq 0$, we get $\Omega = 0$. Otherwise, it has  constant sign. In the last case, we get that $\int_M|\Omega|v_g=2\pi\chi(M),$ which means that the l.h.s. of this equality is a topological invariant.
\hfill$\square$\\ \\
Next, we will give another proof of the B\"ar-type inequality (\ref{barspinc}) for the eigenvalues of any $\Spinc$ Dirac operator. The following theorem was proved by the second author in \cite{r1} using conformal deformation of the spinorial Levi-Civita connection. 
\begin{thm}
 Let $M$ be a compact surface. For any $\Spinc$ structure on $M$, any eigenvalue $\lambda$ of the Dirac operator $D$ to which is attached an eigenspinor $\psi$ satisfies
\begin{eqnarray}
\lambda^2 \geqslant \frac{2\pi \chi(M)}{Area(M)} - \frac{1}{Area(M)} \int_M \vert\Omega\vert v_g.
\label{inequalitybar}
\end{eqnarray}
Equality holds if and only if the eigenspinor $\psi$ is a
 $\Spinc$ Killing spinor, i.e., it satisfies $\Omega\cdot\psi = i \vert\Omega\vert\psi$ and $\nabla_X \psi = -\frac{\lambda}{2} X\cdot\psi$ for any $X \in \Gamma(TM)$.
\label{thmspinc}
\end{thm}
{\bf Proof.} With the help of Remark \eqref{rempropre}, we have that
\begin{eqnarray}
 \lambda^2 = \frac{S}{4} + \vert T^\psi\vert^2 + \bigtriangleup f +(\frac i2 \Omega\cdot\psi, \frac{\psi}{\vert\psi\vert^2}).
\label{barspinccc}
\end{eqnarray}
Substituting the Cauchy-Schwarz inequality, i.e. $\vert T^\psi\vert^2 \geqslant \frac{\lambda^2}{2}$ and the estimate (\ref{cs}) into Equality (\ref{barspinccc}), we easily deduce the result after integrating over $M$. Now the equality in \eqref{inequalitybar} holds if and only if the eigenspinor $\psi$ satisfies $\Omega\cdot\psi = i \vert\Omega\vert\psi$ and $\vert T^\psi\vert^2= \frac{\lambda^2}{2}$. Thus, the second equality is equivalent to say that $\ell^\psi (X) = \frac{\lambda}{2}X$ for all $X\in \Gamma(TM)$.  
Finally, a straightforward computation of the spinorial curvature of the spinor field $\psi$ gives in a local frame $\{e_1,e_2\}$ after using the fact $\beta=-(*\delta)$ that
\begin{eqnarray*}
\frac{1}{2}R_{1212}\ e_1\cdot e_2\cdot\psi&=&\Big(\frac{\lambda^2}{2}+e_1(\delta(e_1))+e_2(\delta(e_2))\Big)e_2\cdot e_1\cdot\psi-\lambda\delta(e_2) e_1\cdot\psi\\&&+\lambda \delta(e_1)e_2\cdot\psi+\Big(e_1(\delta(e_2))-e_2(\delta(e_1))\Big)\psi.
\end{eqnarray*}
Thus the scalar product with $e_1\cdot\psi$ and $e_2\cdot\psi$ implies that $\delta=0$. Finally, $\beta=0$ and the eigenspinor $\psi$ is a $\Spinc$ Killing spinor.
\hfill$\square$\\ \\
Now, we will give some examples where equality holds in (\ref{inequalitybar}) or in \eqref{dett}. Some applications of Theorem \ref{thm1} are also given.\\ \\
{\bf Examples:}
\begin{enumerate}
 \item Let $\cercle^2$ be the round sphere equipped with the standard metric of curvature one. As a K\"ahler manifold, we endow the sphere with the canonical $\Spinc$ structure of curvature form equal to $i\Omega=-i\ltimes$, where $\ltimes$ is the K\"{a}hler $2$-form. Hence, we have $\vert\Omega\vert = \vert\ltimes\vert =1$. Furthermore, we mentionned that for the canonical $\Spinc$ structure, the sphere carries parallel spinors, i.e., an eigenspinor associated with the eigenvalue $0$ of the Dirac operator $D$. Thus equality holds in (\ref{inequalitybar}). On the other hand, the equality in \eqref{dett} also holds, since the sign of the curvature form $\Omega$ is constant. 
\item Let $f: M \rightarrow \cercle^3$ be an isometric immersion of a surface
$M^2$ into the sphere equipped with its unique $\Spin$ structure and assume that the mean curvature $H$ is
constant. The restriction of a Killing spinor on $\cercle^3$ to the surface $M$ defines a spinor field
$\phi$ solution of the following equation \cite{gj}
\begin{eqnarray}
 \nabla_X \phi = -\frac 12 II(X) \bullet \phi + \frac 12 J(X)\bullet\phi,
\label{generalized}
\end{eqnarray}
where $II$ denotes the second fundamental form of the surface and $J$ is the complex
structure of $M$ given by
the rotation of angle $\frac{\pi}{2}$ on $TM$. It is
easy to check that $\phi$ is an eigenspinor for $D^2$ associated with the
eigenvalue $H^2 +1$. Moreover $D \phi = H \phi + e_1\cdot e_2\cdot \phi$, so that $Y=0$. Moreover the tensor $T^\phi=\frac{1}{2}II$ and $Q^\phi=\frac{1}{2}J$. Hence by Theorem \ref{thm1}, and since the norm of $\phi$ is
constant, we obtain
$$H^2 + \frac12 = \frac14 S +\frac{1}{4}|II|^2.$$
\item On two-dimensional manifolds, we can define another Dirac operator associated with the complex
structure $J$ given by $\widetilde D =Je_1 \cdot \nabla_{e_1} +Je_2 \cdot \nabla_{e_2}  = e_2\cdot \nabla_{e_1}- e_1 \cdot\nabla_{e_2}$. Since $\widetilde D$ satisfies $D^2 = (\widetilde D)^2$, all the above results are also true for the eigenvalues of $\widetilde D$.

\item Let $M^2$ be a surface immersed in $\cercle^2 \times \RR$. The product of the canonical $\Spinc$ structure on $\cercle^2$ and the unique $\Spin$ structure on $\RR$ define a $\Spinc$ structure on $\cercle^2\times \RR$ carrying parallel spinors \cite {mo}. Moreover, by the Schr\"{o}dinger-Lichnerowicz formula, any parallel spinor $\psi$ satisfies $\Omega^{\cercle^2\times \mathbb{R}}\cdot\psi=i\psi$, where  $\Omega^{\cercle^2\times \mathbb{R}}$ is the curvature form of the auxiliary line bundle. Let $\nu$ be a unit normal vector field of the surface. We then write $\partial t=T+f\nu$ where $T$ is a vector field on $TM$ with $||T||^2+f^2=1$. On the other hand, the vector field $T$ splits into $T=\nu_1+h\partial t$ where $\nu_1$ is a vector field on the sphere. The scalar product of the first equation by $T$ and the second one by $\partial t$ gives  $||T||^2=h$ which means that $h=1-f^2$. Hence the normal vector field $\nu$ can be written as $\nu=f\partial t-\frac{1}{f}\nu_1.$ As we mentionned before, the $\Spinc$ structure on $\cercle^2 \times \RR$ induces a $\Spinc$ structure on $M$ with induced auxiliary line bundle. Next, we want to prove that the curvature form of the auxiliary line bundle of $M$ is equal to $i\Omega(e_1,e_2)=-if$, where $\{e_1,e_2\}$ denotes a local orthonormal frame on $TM$. Since the spinor $\psi$ is parallel, we have by \cite{mo} that for all $X\in T(\cercle^2 \times \RR)$ the equality ${\rm Ric}^{\cercle^2 \times \RR} X\cdot\psi=i(X\lrcorner\Omega^{\cercle^2 \times \RR})\cdot\psi$. Therefore, we compute
\begin{eqnarray*}
(\nu\lrcorner \Omega^{\cercle^2 \times \RR})\bullet\phi&=&(\nu\lrcorner \Omega^{\cercle^2 \times \RR})\cdot\nu\cdot\psi|_M=i\nu\cdot{\rm Ric}^{\cercle^2 \times \RR}\ \nu\cdot\psi|_M\\
&=&-\frac{1}{f}i\nu\cdot\nu_1\cdot\psi|_M=i\nu\cdot (\nu-f\partial t).\psi|_M\\
&=&(-i\psi-if\nu\cdot\partial t\cdot\psi)|_M. 
\end{eqnarray*}
Hence by Equation \eqref{glucose}, we get that $\Omega\bullet\phi=-i (f\nu\cdot\partial t\cdot\psi)|_M.$ The scalar product of the last equality with $e_1\cdot e_2\cdot\psi$ gives
$$\Omega(e_1,e_2)|\phi|^2=-f\Re(i\nu\cdot\partial t\cdot\psi,e_1\cdot e_2\cdot\psi)|_M= -f\Re(i\partial t\cdot\psi,\psi)|_M.$$
We now compute the term $i\partial t\cdot\psi$. For this, let $\{e'_1, J e'_1\}$ be a local orthonormal frame of the sphere $\cercle^2$. The complex volume form acts as the identity on the spinor bundle of $\cercle^2 \times \RR$, hence $\partial t \cdot\psi = e'_1\cdot Je'_1\cdot\psi$. But we have
$$\Omega^{\cercle^2 \times \RR}\cdot\psi = -\rho\cdot\psi = -\ltimes\cdot\psi = -e'_1\cdot Je'_1\cdot\psi.$$
Therefore, $i\partial t\cdot\psi = \psi$. Thus we get $\Omega(e_1,e_2)=-f$. Finally, 
$$(i\Omega\bullet\phi,\phi) = f \Re(\nu\cdot\partial t\cdot\psi,\psi)|_M=-fg(\nu,\partial t)\vert\phi\vert^2=-f^2\vert\phi\vert^2.$$ Hence Equality in Theorem \ref{thm1} is just 
$$H^2=\frac{S}{4}+\frac{1}{4}|II|^2-\frac{1}{2}f^2.$$
\end{enumerate}

\section{The 3-dimensional case}\label{sec4} 
\setcounter{equation}{0}
In this section, we will treat the $3$-dimensional case.
\begin{thm}
 Let $(M^3,g)$ be an oriented Riemannian manifold. For any $\Spinc$ structure on $M$, any eigenvalue
$\lambda$ of the Dirac operator to which is attached an eigenspinor $\psi$ satisfies
$$\lambda^2 \leqslant \frac{1}{\vol(M, g)} \int_M (\vert T^\psi\vert^2 +\frac S4
+\frac{\vert\Omega\vert}{2})v_g.$$
Equality holds if and only if the norm of $\psi$ is constant and $\Omega\cdot\psi
= i\vert\Omega\vert \psi$. 
\label{thm4}
\end{thm}

{\bf Proof.} As in the proof of Theorem \ref{thm1}, the set $\{\frac{\psi}{|\psi|}, \frac{e_1\cdot\psi}{|\psi|},\frac{e_2\cdot\psi}{|\psi|},\frac{e_3\cdot\psi}{|\psi|}\}$ is orthonormal with respect to the real product $\Re(\cdot, \cdot)$. The covariant derivative of $\psi$ can be expressed in this frame as 
\begin{eqnarray} 
 \nabla_X\psi = \eta(X) \psi + \ell(X)\cdot\psi,
\label{eq1}
\end{eqnarray}
for all vector fields $X,$ where $\eta$ is a $1$-form and $\ell$ is a $(1,1)$-tensor field. Moreover $\eta =\frac{d(\vert\psi\vert^2)}{2\vert\psi\vert^2}$ and $\ell(X)=- \ell^\psi(X)$. Using (\ref{bochner}), it follows that
\begin{eqnarray*}
\lambda^2 &=&  \frac {\Delta ( \vert \psi \vert^2 )} { 2  \vert \psi \vert^2 } +
\vert T^\psi \vert^2 + \frac{\vert d(\vert\psi\vert^2)\vert^2}{4\vert\psi\vert^4}  +\frac 14 S + (\frac i2 \Omega\cdot\psi, \frac{\psi}{\vert\psi\vert^2})\\
&=&   \Delta f -  \frac{\vert d(\vert\psi\vert^2)\vert^2}{2\vert\psi\vert^4} 
+\vert T^\psi \vert^2  +\frac 14 S + (\frac i2 \Omega\cdot\psi, \frac{\psi}{\vert\psi\vert^2}).
 \end{eqnarray*}
By the Cauchy-Schwarz inequality, we have $\frac 12 (i\Omega \cdot\psi,\frac{\psi}{\vert\psi\vert^2}) \leqslant \frac 12 \vert\Omega\vert$. Integrating over $M$ and using the fact that $\vert d(\vert\psi\vert^2)\vert^2\geqslant 0$, we get the result.
\hfill$\square$


\begin{example}
Let $M^3$ be a 3-dimensional Riemannian manifold immersed in $\CC \PP^2$ with
constant mean curvature $H$. Since $\CC \PP^2$ is a K\"{a}hler manifold, we endow it
with the canonical $\Spinc$ structure whose line bundle has
curvature equal to $-3i\ltimes$. Moreover, by the Schr\"{o}dinger-Lichnerowicz formula we have that any parallel spinor $\psi$ satisfies $\Omega^{\CC \PP^2} \cdot\psi = 6i \psi.$ As in the previous example, we compute 
 $$(\nu\lrcorner\Omega^{\CC \PP^2})\bullet\phi = i (\nu\cdot \Ric^{\CC
\PP^2}(\nu)\cdot\psi)_{\mid_M} = -3i \phi.$$ 
Finally, $\Omega\bullet\phi = 3i\phi.$ Using Equation (\ref{diracgauss}), we have that $-\frac 32 H$ is an eigenvalue of $D$. Since the norm of $\phi$ is constant, equality holds in Theorem \ref{thm4} and hence
$$\frac 94 H^2 + \frac 32 = \frac{S}{4}+\frac{1}{4}\vert II\vert^2.$$
\end{example}

\section{Characterization of surfaces in $\mathbb{S}^2\times \mathbb{R}$}\label{charac} 
\setcounter{equation}{0}
In this section, we characterize the surfaces in $\mathbb{S}^2\times \mathbb{R}$ by solutions of the generalized Killing spinors equation which are restrictions of parallel spinors of the canonical $\Spinc$-structure on $\mathbb{S}^2\times \mathbb{R}$ (see also \cite{r} for a different proof). First recall the compatibility equations for characterization of surfaces in $\mathbb{S}^2\times \mathbb{R}$ established by B. Daniel \cite[Thm 3.3]{daniel}: 
\begin{thm} \label{daniel}Let $(M,g)$ be a simply connected Riemannian manifold of dimension $2$, $A:TM\rightarrow TM$ a field of symmetric operator and $T$ a vector field on $TM$. We denote by $f$ a real valued function such that $f^2+||T||^2=1$. Assume that $M$ satisfies the Gauss-Codazzi equations, i.e. $G={\rm det} A+f^2$ and
$$  d^\nabla A(X,Y):=(\nabla_XA)Y-(\nabla_YA)X=f(g(Y,T)X-g(X,T)Y),$$ 
where $G$ is the gaussian curvature, and the following equations 
$$\nabla_X T=fA(X), \,\,\, X(f)=-g(AX,T).$$ Then there exists an isometric immersion of $M$ into $\mathbb{S}^2\times \mathbb{R}$ such that the Weingarten operator is $A$ and  $\partial t=T+f\nu,$ where $\nu$ is the normal vector field to the surface $M$. 
\end{thm}
Now using this characterization theorem, we state our result:
\begin{thm} \label{carac}
Let $M$ be an oriented  simply connected Riemannian manifold of dimension $2$. Let $T$ be a vector field and denote by $f$ a real valued function such that $f^2+||T||^2=1$. Denote by $A$ a symmetric endomorphism field of $TM$. The following statements are equivalent: 
\begin{enumerate}
\item There exists an isometric immersion of $M$ into $\mathbb{S}^2\times \mathbb{R}$ of Weingarten operator $A$ such that $\partial t=T+f\nu,$ where $\nu$ is the unit normal vector field of the surface.
\item There exists a $\Spinc$ structure on $M$ whose line bundle has a connection of curvature given by $i\Omega=-if \ltimes,$ such that it carries a non-trivial solution $\phi$ of the generalized Killing spinor equation $\nabla_X\phi=-\frac{1}{2}AX\bullet\phi$, with $T\bullet\phi=-f\phi+\bar\phi.$
\end{enumerate}
\end{thm}
{\bf Proof.} We begin with $1\Rightarrow 2$. The existence of such a $\Spinc$ structure is assured by the restriction of the canonical one on $\cercle^2\times \mathbb{R}$. Moreover, using the spinorial Gauss formula \eqref{spingauss}, any parallel spinor $\psi$ on $\cercle^2\times \mathbb{R}$ induces a generalized Killing spinor $\phi=\psi|_M$ with $A$ the Weingarten map of the surface $M$. Hence it remains to show the relation $T\bullet\phi=-f\phi+\bar\phi$. In fact, using that $\Omega^{\cercle^2\times \mathbb{R}}\cdot\psi=i\psi,$ we write in the frame $\{e_1,e_2,\nu\}$ 
\begin{equation}\label{eq:3}
\Omega^{\cercle^2\times \mathbb{R}}(e_1,e_2)e_1\cdot e_2\cdot\psi+\Omega^{\cercle^2\times \mathbb{R}}(e_1,\nu)e_1\cdot\nu\cdot\psi+\Omega^{\cercle^2\times \mathbb{R}}(e_2,\nu)e_2\cdot\nu\cdot\psi=i\psi.
\end{equation}
By the previous example in Section \ref{sect:3}, we know that $\Omega^{\cercle^2\times \mathbb{R}}(e_1,e_2)=-f$. For the other terms, we compute 
$$
\Omega^{\cercle^2\times \mathbb{R}}(e_1,\nu)=\Omega^{\cercle^2\times \mathbb{R}}(e_1,\frac{1}{f}\partial t-\frac{1}{f}T)=-\frac{1}{f}g(T,e_2)\Omega^{\cercle^2\times \mathbb{R}}(e_1,e_2)=g(T,e_2),
$$
where the term $\Omega^{\cercle^2\times \mathbb{R}}(e_1,\partial t)$ vanishes since we can split $e_1$ into a sum of vectors on the sphere and on $\mathbb{R}.$ Similarly, we find that $\Omega^{\cercle^2\times \mathbb{R}}(e_2,\nu)=-g(T,e_1).$ By substituting these values into \eqref{eq:3} and taking Clifford multiplication with $e_1\cdot e_2$, we get the desired property. For $2\Rightarrow 1$, a straightforward computation for the spinorial curvature of the generalized Killing spinor $\phi$ yields on a local frame $\{e_1,e_2\}$ of $TM$ that 
\begin{equation}
(-G+{\rm det}\,A)e_1\bullet e_2\bullet\phi=-(d^\nabla A)(e_1,e_2)\bullet \phi + if\phi.
\label{eq:6}
\end{equation}
In the following, we will prove that the spinor field $\theta:=i\phi-if\bar{\phi}+JT\bullet\phi$ is zero. For this, it is sufficient to prove that its norm vanishes. Indeed, we compute 
\begin{equation}
|\theta|^2=|\phi|^2+f^2|\phi|^2+||T||^2|\phi|^2-2\Re(i\phi,if\bar{\phi})+2\Re(i\phi,JT\bullet\phi)
\label{eq:5}
\end{equation}
From the relation $T\bullet\phi=-f\phi+\bar{\phi}$ we deduce that $\Re(\phi,\bar\phi)=f|\phi|^2$ and the equalities
$$g(T,e_1)|\phi|^2=\Re(ie_2\bullet\phi,\phi) \quad\text{and}\quad g(T,e_2)|\phi|^2=-\Re(ie_1\bullet\phi,\phi).$$
Therefore, Equation \eqref{eq:5} becomes 
\begin{eqnarray*}
|\theta|^2&=&2|\phi|^2-2f^2|\phi|^2+2\Re(i\phi,JT\bullet\phi)\\
&=&2|\phi|^2-2f^2|\phi|^2+2 g(JT,e_1)\Re(i\phi,e_1\bullet\phi)+2g(JT,e_2) \Re(i\phi,e_2\bullet\phi)\\
&=&2|\phi|^2-2f^2|\phi|^2+2 g(JT,e_1) g(T,e_2)|\phi|^2-2g(JT,e_2) g(T,e_1)|\phi|^2\\
&=&2|\phi|^2-2f^2|\phi|^2-2g(JT,e_1)^2|\phi|^2-2g(T,e_1)^2|\phi|^2\\
&=&2|\phi|^2-2f^2|\phi|^2-2||T||^2|\phi|^2=0.
\end{eqnarray*}
Thus, we deduce $if\varphi=-f^2e_1\cdot e_2\cdot\varphi-fJT\cdot\varphi$, where we use the fact that $\bar\phi=i e_1\bullet e_2\bullet\phi$. In this case, Equation \eqref{eq:6} can be written as 
$$(-G+{\rm det}\,A+f^2)e_1\bullet e_2\bullet\phi=-((d^\nabla A)(e_1,e_2)+fJT)\bullet \phi.$$
This is equivalent to say that both terms $R_{1212}+{\rm det}\,A+f^2$ and $(d^\nabla A)(e_1,e_2)+fJT$ are equal to zero. In fact, these are the Gauss-Codazzi equations in Theorem \ref{daniel}. In order to obtain the two other equations, we simply compute the derivative of $T\cdot\varphi=-f\varphi+\bar\varphi$ in the direction of $X$ to get 
\begin{eqnarray*}
\nabla_X T \bullet \phi+T\bullet\nabla_X\varphi&=&\nabla_X T \bullet \phi-\frac{1}{2}T\bullet A(X)\bullet\varphi\\
&=&-X(f)\varphi-f\nabla_X\varphi+\nabla_X\bar\varphi\\
&=&-X(f)\varphi+\frac{1}{2}f AX\bullet\varphi+\frac{1}{2} AX\bullet\bar{\varphi}\\
&=&-X(f)\varphi+\frac{1}{2}f AX\bullet\varphi+\frac{1}{2} AX\bullet(T\bullet\varphi+f\varphi).
\end{eqnarray*}
This reduces to $\nabla_X T \bullet \phi+g(T,A(X))\varphi=-X(f)\varphi+fA(X)\bullet\varphi.$ Hence we obtain $X(f)=-g(A(X),T)$ and $\nabla_X T=fA(X)$ which finishes the proof.
\hfill$\square$
\begin{remark} The second condition in Theorem \ref{carac} is equivalent to the existence of a $\Spinc$ structure whose line bundle $L$ verifies 
$c_1(L)=[\frac{i}{2\pi}f\ltimes]$ and $f \ltimes$ is a closed 2-form. This $\Spinc$ structure carries a non-trivial solution $\phi$ of the generalized Killing spinor equation $\nabla_X\phi=-\frac{1}{2}AX\bullet\phi$, with $T\bullet\phi=-f\phi+\bar\phi.$
\end{remark}
{\bf Acknowledgment}\\\\
Both authors are grateful to Oussama Hijazi for his encouragements and relevant remarks.

\end{document}